\documentstyle[12pt]{article}
\newtheorem{theor}{Theorem}
\newtheorem{lem}{Lemma}

\newtheorem{prop}{Proposition}
\begin{document}

\title{Properties of parallelotopes equivalent to Voronoi's conjecture}

\author{
Michel Deza, Ecole Normale Sup\'erieure, Paris\\
Viacheslav Grishukhin, CEMI RAN, Moscow, Russia}

\date{}

\maketitle

\begin{abstract}
A parallelotope is a polytope whose translation copies fill space
without gaps and intersections by interior points. Voronoi conjectured
that each parallelotope is an affine image of the Dirichlet domain
of a lattice, which is a Voronoi polytope. We give several properties
of a parallelotope and prove that each of them is equivalent to it is
an affine image of a Voronoi polytope.
\end{abstract}

\section{Introduction}
Let us have a partition of ${\bf R}^n$ into equal convex polytopes
({\em tiles})
such that the partition is simultaneously a covering and a packing and
the intersection of any two polytopes is empty or a common face of each.
Such a partition is called {\em tiling}. Let this tiling be invariant under
a group $\cal T$ of translations, and the group $\cal T$ is transitive on
polytopes of the tiling. Then each tile of such a partition is called
{\em parallelotope}. Here the prefix {\em parallelo} emphasizes that each
tile is a parallel translation of a prototile. (Following to \cite{MRS}
we use the word {\em parallelotope} rather than {\em parallelohedron} which
was used by Voronoi in \cite{Vo}. Recall also that a {\em polyhedron} is a
3-dimensional {\em polytope}.)

G.Voronoi in \S8 of \cite{Vo} defines a parallelohedron as follows. A
polytope $P$ with a group of translations $\cal T$ is called a
parallelohedron if the space ${\bf R}^n$ can be filled by non-overlapping
congruent copies of $P$ using translations taken from $\cal T$.

Each parallelotope necessarily satisfies the following three conditions:

(i) a parallelotope is centrally symmetric;

(ii) each facet of a parallelotope is centrally symmetric;

(iii) for $n>1$, the projection of a parallelotope along any
$(n-2)$-dimensional face is either a parallelogram or a centrally
symmetric hexagon.

The edges of the polygon of item (iii) above are projections of four
or six facets of the projected parallelotope $P$. These facets form a
{\em belt} of the parallelotope $P$. Hence the property (iii) of
a parallelotope $P$ has another formulation. Namely,

(iii)' for $n>1$, each belt of a parallelotope contains four or six facets.

B.Venkov \cite{Ve} (and independently P.McMullen \cite{McM}) proved
that the above three conditions are sufficient for a polytope to be a
parallelotope. A.Aleksandrov \cite{Al}, knowing the Venkov's result,
simplified the proof of B.Venkov.

There is a special well known case of a parallelotope, namely, the Voronoi
polytope related to a point of a lattice $L$. The Voronoi polytope
$P_V(t_0)$ related to a point $t_0\in L$ is the set of points of
${\bf R}^n$ which are at least as close to $t_0$ as to any other point
of $L$.

The main conjecture of Voronoi is that any parallelotope $P$ can be mapped
into a Voronoi polytope $P_V$ under an affine transformation $x\to Ax$ of
the space containing $P$. Here $A$ is a non-degenerate $n\times n$ matrix,
where $n$ is dimension of the space.

Call a $k$-face (=$k$-dimensional face) of a parallelotope {\em primitive}
if it belongs to minimal possible number $n-k+1$ of parallelotopes of its
tiling. Obviously any facet of a parallelotope is primitive.
According to Zhitomorskii \cite{Zh}, a parallelotope is called
{\em k-primitive} if each its $k$-face is primitive. Besides,
the $k$-primitivity implies the $(k+1)$-primitivity. A 0-primitive
parallelotope is simply called {\em primitive}. Obviously, any
parallelotope is $(n-1)$-primitive.

Voronoi proved his conjecture for primitive parallelotopes. If a
parallelotope is primitive, then each its belt
contains 6 facets, but not conversely. On the other hand, if each belt
consists of 6 facets, then the parallelotope is $(n-2)$-primitive. This
implies that each $(n-2)$-face belongs to three parallelotopes.
Zhitomirskii \cite{Zh} extend the result of Voronoi over $(n-2)$-primitive
parallelotopes.

P.McMullen \cite{McM1} proved Voronoi's conjecture for parallelotopes
which are zo\-no\-to\-pes. Later R.Erdahl \cite{Er} repeated the result of
McMullen. B.Delaunay \cite{De} proved Voronoi's conjecture in complete
generality for the dimensions $n\le 4$.

We give here several conditions on a parallelotope $P$ each of which is
equivalent to the Voronoi's conjecture is true for $P$.

\section{Parallelotopes}
Now we consider a description of a parallelotope $P=P(0)$ with the
center in origin by a system of linear inequalities. We denote by $A^T$
the transpose of a matrix $A$ and by $x^Ty=y^Tx$ the scalar product of
two column vectors $x,y\in {\bf R}^n$, and set $x^Tx=x^2$.

Being a convex polytope a parallelotope is described by a system of linear
inequalities $\{q_i^Tx\le \alpha_i\}$. Since, by (i) of Introduction, $P$
is centrally symmetric, each facet $F_i$ of $P$ has the opposite facet
$-F_i$. If $F_i$ lies in the affine hyperplane given by the equality
$q_i^Tx=\alpha_i$, then the opposite facet $-F_i$ lies in an affine
hyperplane $q_i^Tx=-\alpha_i$. The parallelotope $P$ is adjacent by the
facet $F_i$ to a parallelotope $P(t_i)$ which is a parallel shift
of $P$ by the translation vector $t_i\in \cal T$.

Let ${\cal I}_P$ be the set of indices of pairs of opposite facets of $P$.
Then the set $\{t_i:i\in {\cal I}_P\}$ of translation vectors generates the
translation group $\cal T$ and a lattice $L$. The points of $L$, i.e. the
centers of parallelotopes of the tiling, are the end-points of lattice
vectors. We can identify each lattice vector with an element $t\in \cal T$.
By this identification origin is the zero point and simultaneously the
zero lattice vector 0 of $L$.

Obviously, the point $\frac{1}{2}t_i$ is the center of the facet $F_i$.
Hence $\frac{1}{2}q_i^Tt_i=\alpha_i$. We have
\begin{equation}
\label{par}
P(0)=\{x\in{\bf R}^n:-\frac{1}{2}q_i^Tt_i\le q_i^Tx\le\frac{1}{2}q_i^Tt_i,
\mbox{  }i\in{\cal I}_P\}.
\end{equation}

Here the {\em facet vectors} $q_i$ are determined only
up to a non-zero multiple $\beta_i$. We say that the facet vectors $q_i$
and the lattice vectors $t_i$, $i\in{\cal I}_P$, giving the description
(\ref{par}) of a parallelotope $P$, are {\em associated}.

The parallelotope $P(t)$ with the center in the point $t\in L$ is a
translation of $P(0)$ by the vector $t$, and therefore it is described 
as follows:
\begin{equation}
\label{pat}
P(t)=\{x\in{\bf R}^n:-\frac{1}{2}q_i^Tt_i\le q_i^T(x-t)
\le\frac{1}{2}q_i^Tt_i, \mbox{  }i\in{\cal I}_P\}.
\end{equation}

Note that, by definition of a Voronoi polytope, a facet $F_i$ of the
Voronoi polytope $P_V(0)$ is orthogonal to the lattice vector $t_i$ and
bisects it. The lattice vector $t_i$ connects the centers of $P_V(0)$ and
$P_V(t_i)$, where $P_V(t_i)$ is the Voronoi polytope adjacent to $P_V(0)$
by the facet $F_i$. In other words, we can set $q_i=t_i$ in the
descriptions (\ref{par}) and (\ref{pat}) of parallelotopes $P(0)$ and
$P(t)$ in the case they are Voronoi polytopes.

\section{Linear transforms of parallelotopes}
Note that the usual Euclidean norm $x^2$ is used in the definition
of the Voronoi polytope $P_V(0)$. But we can use an arbitrary positive
quadratic form $f(x)=x^TDx$ as a norm of $x$. Here $D$ is a symmetric
positive definite $n\times n$ matrix. Then the above definition
gives a parallelotope $P_f$. Call such a parallelotope the {\em Voronoi
polytope with respect to the quadratic form $f(x)$}. Such a parallelotope
relates to a lattice $L$ (or to a translation group $\cal T$). Consider the
Voronoi polytope $P_f$ with respect to the quadratic form $x^TDx$ in detail.
By definition, we have
\[P_f(0)=\{x\in{\bf R}^n:x^TDx\le(x-t)^TD(x-t)\mbox{ for all }
t\in {\cal T}\}.\]
It is well known, that a finite set $\{\pm t_i:i\in {\cal I}_f\}$ of
vectors $t_i\in \cal T$ is sufficient for the description of $P_f(0)$.

Using the symmetricity of $D$ and joining the inequalities for $t_i$ and
$-t_i$, we simplify this as follows
\begin{equation}
\label{vor}
P_f(0)=\{x\in {\bf R}^n:-\frac{1}{2}t_i^TDt_i\le t_i^TDx\le
\frac{1}{2}t_i^TDt_i, \mbox{  }i\in {\cal I}_f\}.
\end{equation}
For $D=I_n$, where $I_n$ is the identity matrix, we have $f(x)=x^2$ and
$P_f(0)=P_V(0)$.
\begin{lem}
\label{ADq}
Let $P$ be a parallelotope given by (\ref{par}). Let $A$ be an $n\times n$
non-degenerate matrix, and $D=A^TA$. The following assertions are
equivalent:

(i) the affine transformation $x\to Ax$ transforms $P$ into a Voronoi
polytope;

(ii) $P$ is the Voronoi polytope with respect to the quadratic form
$f(x)=x^TDx$;

(iii) the facet vectors $q_i$ satisfy the equality $q_i=Dt_i$,
$i\in{\cal I}_P$.
\end{lem}
{\bf Proof}. (i)$\Rightarrow$(iii).
Consider the affine transformation $x\to Ax$. 
The new facet vector has the form $(A^T)^{-1}q$. In fact, for $x\in F$,
the point $Ax$ belongs to a facet of $AP$. Hence
\[((A^T)^{-1}q)^TAx=q^T((A^T)^{-1})^TAx=q^TA^{-1}Ax=q^Tx=0. \]
The new center of the transformed facet is $\frac{1}{2}At_i$. For $AP$
to be a Voronoi polytope, we have to have
\[(A^T)^{-1}q_i=At_i, \mbox{ i.e. }q_i=A^TAt_i, \mbox{ i.e. }q_i=Dt_i. \]

(ii)$\Rightarrow$(i)
The positive definite matrix $D$ can be represented as the product
$D=A^TA$, where the matrix $A$ is non-degenerate. Hence the form
$x^TDx=x^TA^TAx=(Ax)^T(Ax)$ is the quadratic form $(Ax)^2$ in the
transformed space, i.e. $P=P_f=P_{(Ax)^2}$. Let $y=Ax$. Then $x=A^{-1}y$
and $P=P_{(Ax)^2}=A^{-1}P_{y^2}=A^{-1}P_V$. Hence $AP=P_V$.

(iii)$\Rightarrow$(ii)
If we set in (\ref{par}) $q_i=Dt_i$, we obtain description (\ref{vor}) of
a parallelotope. Hence the Voronoi polytope with respect to a quadratic
form is also a special case of a parallelotope, when $q_i=Dt_i$.
In other words, in this case, the parallelotope $P(0)$ of (\ref{par})
is $P_f(0)$ for $f(x)=x^TDx$. \hfill $\Box$

\section{A canonical representation of a parallelotope}

Consider a vertex $v$ of a facet $F_i$. Let $v$ be the intersection
of facets $F_j$, $j\in{\cal I}(v)$. Then $i\in{\cal I}(v)$ and
\[q_j^Tv=\frac{1}{2}q_j^Tt_j, \mbox{  }j\in{\cal I}(v). \]
Since each facet of $P$ is centrally symmetric, there is a symmetric to
$v$ vertex $v^s\in F_i$. We have
\[q_k^Tv^s=\frac{1}{2}q_k^Tt_k, \mbox{  }k\in{\cal I}(v^s). \]
Note that the point $\frac{1}{2}(v+v^s)$ is the center $\frac{1}{2}t_i$
of the facet $F_i$. Hence
\[v+v^s=t_i, \mbox{ i.e. }v^s=t_i-v. \]

Recall that there are the following two types of belts in the
parallelotope $P$:

1) 3-belts containing 6 facets $\pm F_i$, $\pm F_j$, $\pm F_k$;

2) 2-belts containing 4 facets $\pm F_i$, $\pm F_j$.

We denote each belt by the set of indices of its generating facets. So,
we have the following two types of belts: $\{i,j,k\}$ and $\{i,j\}$.

Therefore some facet vectors of the set $\{q_j:j\in{\cal I}(v)\}$ and
$\{q_k:k\in{\cal I}(v^s)\}$ are joined into pairs of two types such that

1) either the facets $F_i,F_j,F_k$ belong to the 3-belt $\{i,j,k\}$;

2) or the facets $F_j$ and $F_k$ are opposite, i.e. $F_k=-F_j$, and
belong to the 2-belt $\{i,j\}$.

Let ${\cal B}$ be the set of all 3-belts. Consider a belt
$\{i,j,k\}\in {\cal B}$. The facets vectors
$q_i, q_j, q_k$ lie in a two-dimensional plane, where they are linearly
dependent. Obviously, the associated lattice vectors $t_i,t_j,t_k$ are
also linearly dependent. Moreover, this dependence has the following form
\begin{equation}
\label{et}
t_i-\varepsilon_jt_j-\varepsilon_kt_k=0, \mbox{  }\{i,j,k\}\in{\cal B},
\end{equation}
where $\varepsilon_j,\varepsilon_k\in\{\pm 1\}$. Since each facet vector
$q_i$ is defined up to a scalar multiple $\beta_i$, we can choose lengths
of the associated facet vectors such that the new facet vector
$\beta_iq_i$ satisfy the equality similar to (\ref{et})
\begin{equation}
\label{eq}
\beta_iq_i-\varepsilon_j\beta_jq_j-\varepsilon_k\beta_kq_k=0,
\mbox{  }\{i,j,k\}\in{\cal B}.
\end{equation}
Following to \cite{Vo} and \cite{RR}, we say that, for the belt $\{i,j,k\}$,
the facet vectors $q_i,q_j,q_k$ are defined {\em canonically with respect
to the 3-belt} $\{i,j,k\}$ if they satisfy the same equality as the
associated lattice vectors $t_i,t_j,t_k$.

\vspace{3mm}
{\bf Definition} A parallelotope $P$ {\em is defined canonically} by
(\ref{par}) if the facet vectors $q_i$, $i\in{\cal I}_P$, are defined
canonically simultaneously with respect to all belts of $P$.

In other words, a parallelotope $P$ {\em is defined canonically} if the
system of equations (\ref{eq}) determining multiples $\beta_i$,
$i\in{\cal I}_P$, has a non-zero solution.

\vspace{2mm}
Voronoi proves in the first \S\S 44 of \cite{Vo} that a primitive
parallelotope can be defined canonically.

\section{Relations between lattice and facet vectors}

We suppose that the facet vectors $q_i$, $q_j$ and $q_k$ determine the
facets $F_i$, $F_j$ and $F_k$, respectively. Hence the vector
$\varepsilon_jq_j$ defines the facet $\varepsilon_jF_j$.

In Proposition~\ref{belt} below, for simplicity sake, we suppose that
$\varepsilon_j=\varepsilon_k=1$. For to apply the results below to the
general dependences (\ref{et}) and (\ref{eq}), it is sufficient to
change $q_j$ and $q_k$ by $\varepsilon_jq_j$ and $\varepsilon_kq_k$,
respectively.

So the lattice vecors, corresponding to the belt $\{i,j,k\}$,
satisfy the equality
\begin{equation}
\label{tijk}
t_i=t_j+t_k.
\end{equation}
Hence the defined canonically facet vectors satisfy the equality
\begin{equation}
\label{qijk}
q_i=q_j+q_k
\end{equation}
and the intersections $F_i\cap F_j$ and $F_i\cap F_k$ are non-empty
and define two opposite facets of $F_i$.

For $i\in{\cal I}_P$, let
${\cal I}_i=\{j\in{\cal I}_P:F_j\cap F_i=F^{n-2}\}$, where $F^{n-2}$ is
an $(n-2)$-face of $P$ (which is a facet of $F_i$). Let ${\cal B}_i$ be
the set of 3-belts containing the facet $F_i$.
Now, using the property (ii) of parallelotopes, we prove an important fact.
\begin{prop}
\label{belt}
For all $j\in{\cal I}_i$, the vectors $q_j$ can be defined canonically
with respect to all 3-belts of ${\cal B}_i$. For these canonical facet
vectors we have

1) for a 3-belt $\{i,j,k\}\in{\cal B}_i$ the following equalities hold:
\[q_i^Tt_j=q_j^Tt_i, \mbox{  }q_i^Tt_k=q_k^Tt_i, \mbox{  }
q_j^Tt_k=q_k^Tt_j; \]

2) for a 2-belt $\{i,j\}$, $j\in{\cal I}_i$ the following equalities hold:
\[q_j^Tt_i=q_i^Tt_j=0. \]
\end{prop}
{\bf Proof}. Obviously, the vectors $q_j$, $j\in{\cal I}_i$, can be
defined canonically with respect to all 3-belts of ${\cal B}_i$, since
the belts of ${\cal B}_i$ have only one common vector $q_i$.

Recall that opposite vertices $v$ and $v^s$ of the facet $F_i$ are
determined by facet vectors some of which form belts with $q_i$.
Consider the equations
\[q_j^Tv=\frac{1}{2}q_j^Tt_j, \mbox{  }q_k^Tv^s=\frac{1}{2}q_k^Tt_k, \]
related either to a 3-belt $\{i,j,k\}$ or to a 2-belt $\{i,j\}$, and
then $q_k=-q_j$.

For the case 1), substituting the expressions $q_k=q_i-q_j$,
$t_k=t_i-t_j$, $v^s=t_i-v$, in the second equation, and using the
equalities $q_i^Tv=\frac{1}{2}q_i^Tt_i$, $q_j^Tv=\frac{1}{2}q_j^Tt_j$,
we obtain
\[q_k^Tv^s=\frac{1}{2}q_k^Tt_k \Rightarrow
(q_i-q_j)^T(t_i-v)=\frac{1}{2}(q_i-q_j)^T(t_i-t_j)\Rightarrow \]
\[q_i^T(t_i-v-\frac{1}{2}t_i+\frac{1}{2}t_j)=
q_j^T(t_i-v-\frac{1}{2}t_i+\frac{1}{2}t_j)\Rightarrow
q_i^Tt_j=q_j^Tt_i. \]
Similarly, beginning with $q_j^Tv=\frac{1}{2}q_j^Tt_j$ and using the
equality $v=t_i-v^s$, we obtain the equality $q_i^Tt_k=q_k^Tt_i$. Now,
this equality implies
\[(q_j+q_k)^Tt_k=q_k^T(t_j+t_k)\Rightarrow q_j^Tt_k=q_k^Tt_j. \]

In the case 2) we have $q_k=-q_j$, $t_k=-t_j$. Hence we obtain
\[q_k^Tv^s=\frac{1}{2}q_k^Tt_k\Rightarrow
q_j^T(v-t_i)=\frac{1}{2}q_j^Tt_j\Rightarrow
 q_j^T(\frac{1}{2}t_j-t_i)=\frac{1}{2}q_j^Tt_j\Rightarrow
q_j^Tt_i=0. \]
Using the facet $F_j$ instead of $F_i$, we obtain the equality $q_i^Tt_j=0$.
Note that the equalities $q_j^Tt_i=0=q_i^Tt_j$ do not depend on whether
$q_i$ and $q_j$ are canonical or not. \hfill $\Box$

Let $|{\cal I}_P|=m$ and let $Q$ and $T$ be $n\times m$ matrices whose
columns are the vectors $q_i$ and $t_i$ for $i\in{\cal I}_P$, respectively.
Then the product $q_i^Tt_j$ is the $(ij)$-th element of the matrix
product $Q^TT$. If the equalities
\begin{equation}
\label{QT}
q_i^Tt_j=t_i^Tq_j \mbox{ hold for all pairs }i,j\in {\cal I}_P,
\end{equation}
then the $m\times m$ matrix $Q^TT$ is symmetric, i.e. $Q^TT=(Q^TT)^T=T^TQ$.

\begin{lem}
\label{qt}
The following assertions are equivalent:

(i) the equalities $q_i^Tt_j=t_i^Tq_j$ hold for all pairs
$i,j\in{\cal I}_P$;

(ii) there is a unique symmetric non-degenerate $n\times n$ matrix
$D$ such that $q_i=Dt_i$ for all $i\in{\cal I}_P$.
\end{lem}
{\bf Proof}. (i) $\Rightarrow$ (ii) Let ${\cal I}_b\subseteq{\cal I}_P$
be an $n$-subset of ${\cal I}_P$ such that the set $\{t_i:i\in{\cal I}_b\}$
is linearly independent. Let $T_b$ and $Q_b$ be the submatrices of $T$ and
$Q$ composed by column vectors $t_i$ and $q_i$ for $i\in{\cal I}_b$,
respectively. Then $T_b$ is an $n\times n$ nondegenerate matrix. If
(\ref{QT}) is true, then it implies the equality $T_b^TQ=Q_b^TT$. This
equality is equivalent to the equality
\[Q=DT, \]
where $D=(T_b^T)^{-1}Q_b^T=(Q_bT_b^{-1})^T$. The matrix $D$ is
symmetric. In fact, take a restriction of the equality $Q=DT$ onto
the columns $t_i,q_i$ for $i\in {\cal I}_b$. The restriction is
$Q_b=DT_b$, i.e. $Q_bT_b^{-1}=D=(Q_bT_b^{-1})^T$. The matrix $D$ does not
depend on a chose of $T_b$. In fact, if there is another symmetric matrix
$D'$ such that $Q=D'T$, then $D=Q_bT_b^{-1}=D'T_bT_b^{-1}=D'$.

(ii) $\Rightarrow$ (i) Conversely, let $q_i=Dt_i$. Then
$q_i^Tt_j=(Dt_i)^Tt_j=t_i^TDt_j=t_i^Tq_j$. \hfill $\Box$

\section{Graphs related to tilings}

Recall that the centers of parallelotopes $P(t)$ form a lattice $L$.
Consider the points of $L$ (i.e. the endpoints of lattice vectors) as
vertices of a graph $G_L$. Two vertices $t,t'\in L$ are adjacent in $G_L$
if and only if $t-t'\in\{\pm t_i:i\in {\cal I}_P\}$. We can consider
$G_L$ as an oriented graph, where the orientation of the edge $t_i$ is
the orientation of the vector $t_i$. Hence edges of $G_L$
are the vectors $\pm t_i$, $i\in {\cal I}_P$. Therefore the set of all
edges of $G_L$ is partitioned into $m=|{\cal I}_P|$ classes $E_i$. We
suppose that all edges of $E_i$ are vectors $t_i$ (with the same sign,
say +). In other words, all edges of $E_i$ are obtained from one by
translations. For each $i\in{\cal I}_P$, a vertex of $G_L$ is incident
to two edges $t_i$, one of which go out and another come in the vertex.

For any collection of integers $\{z_i:i\in{\cal I}_P\}$, we set
$t(z)=\sum_{i\in{\cal I}_P}z_it_i$ and
$q(z)=\sum_{i\in{\cal I}_P}z_iq_i$. Let
${\cal I}(z)=\{i\in{\cal I}_P:z_i\not=0\}$ be the {\em support} of 
$z=\{z_i:i\in{\cal I}_P\}$.

Any two vertices $t^0$ and $t$ are connected in $G_L$ by an oriented path
$\cal P$ directed from $t^0$ to $t$. We denote such a path as a sequence
${\cal P}=(t_1,t_2,...,t_s)$ of vectors $t_k\in\{\pm t_i:i\in{\cal I}_P\}$
corresponding to edges of the path in the natural order. Here $t_k=t_i$ or
$t_k=-t_i$ according to the directions of the path and the corresponding
edge coincide or not, respectively. Then, obviously, $t=t^0+\sum_{k=1}^st_k$.
In particular, if the path is closed, i.e. it is a circuit and $t=t^0$,
then $\sum_{k=1}^st_k=0$. We rewrite the sum $\sum_{k=1}^st_k$ as
$\sum_{i\in{\cal I}_P}z_i({\cal P})t_i=t(z({\cal P}))=t({\cal P})$. So,
$t=t^0+t({\cal P})$. We set ${\cal I}({\cal P})={\cal I}(z({\cal P}))$.
Note that there are many paths with the same collection $z({\cal P})$.
All of them are obtained from $\cal P$ by permutations the edges $t_k$.

Since the set $\{\pm t_i:i\in {\cal I}_P\}$ generates the lattice $L$,
any lattice vector $t$ has a (non-unique) representation $t=t({\cal P})$,
where $\cal P$ is a path in $G_L$ connecting 0 with $t$. We associate
the vector $q({\cal P})=\sum_{t_k\in{\cal P}}q_k$ to the vector
$t({\cal P})$. We call the vector $q({\cal P})$ the
{\em vector associated} to the vector $t({\cal P})$.
If the equalities $q_i^Tt_j=q_j^Tt_i$ hold for all $i,j\in {\cal I}_P$,
then using Lemma~\ref{qt} we see that $q({\cal P})=Dt({\cal P})$ does not
depend on $\cal P$. The uniqueness of $q$ in this case follows from the
fact that the equality $\sum_kt_k=0$ implies the equality $\sum_kq_k=0$.

Consider some important subgraphs of the graph $G_L$. Let $t\in L$ and
let $G(t)$ be the graph induced by all vertices $t'\in L$ adjacent to $t$.
The union of $t$ and $G(t)$ is the suspension $\nabla G(t)$. In the graph
$\nabla G(t)$, the vertex $t$ is adjacent to all vertices of $G(t)$.

For $t=0$, the vertices of $G(0)$ are endpoints of the vectors
$\pm t_i$, $i\in{\cal I}_P$. In other words, the graph $G(0)$ is
determined on two copies of the set ${\cal I}_P$. In fact, any edge
$t_i$ of $G(0)$ belongs to a triangle $\triangle=(t_i,t_j,t_k)$ of
$\nabla G(0)$ such that $t_i-\varepsilon_jt_j-\varepsilon_kt_k=0$.
This triangle corresponds to the belt $\{i,j,k\}$ of $P(0)$.
Moreover, there are six such triangles in $\nabla G(0)$. These 6
triangles form a hexagon spanning a 2-dimensional plane.

Let $G(F^k)$ be the graph induced by the centers of all parallelotopes
having a common $k$-face $F^k$ with $P(0)$. So, the graph $G(F^{n-1})$
is an edge $t_i$ with end-vertices corresponding to adjacent
parallelotopes $P(0)$ and $P(t_i)$. There are only two types of graphs
$G(F^{n-2})$, namely, triangles and quadrangles, according to two types
of belts with 6 and 4 facets, respectively.

If $P$ is $m$-primitive, then, for $m\le k\le n$, $G(F^k)=K_{n-k+1}$,
where $K_s$ is the complete graph on $s$ vertices.

The following reformulation of canonical definity of a parallelotope $P$
is obvious.
\begin{lem}
\label{tri}
For a parallelotope $P$, the following assertions are equivalent:

(i) $P$ is defined canonically;

(ii) $q(\triangle)=0$ for all triangles $\triangle\subset G_L$ such that
$\triangle=G(F^{n-2})$.   \hfill $\Box$
\end{lem}

Call a 4-circuit $(t_i,t_j,-t_i,-t_j)$ by a {\em quadrangle}
${\cal Q}_{ij}$. It is a parallelogram and spans a 2-dimensional plane. 
Among quadrangles of $G_L$ there are quadrangles $G(F^{n-2})$, where 
$F^{n-2}$ is an $(n-2)$-face which is common to four parallelotopes.
Obviously, for each quadrangle, we have trivially $q({\cal Q}_{ij})=
q_i+q_j-q_i-q_j=0$.

The technique used in \cite{Ve} can be applied for a proof of
Proposition~\ref{Ven} below (see also Theorem 1 of \cite{RR}).
\begin{prop}
\label{Ven}
Any circuit of $G_L$ can be represented as a sum modulo 2 of circuits
of type $G(F^{n-2})$.
\end{prop}
Hence Lemma~\ref{tri} and Proposition~\ref{Ven} imply the following
\begin{lem}
\label{cac}
For a parallelotope $P$, the following assertions are equivalent:

(i) $P$ is defined canonically;

(ii) $q({\cal C})=0$ for all circuits ${\cal C}\subset G_L$.
\end{lem}

But we give here an explicit proof of a more weak result which we'll
use later.
\begin{lem}
\label{quad}
Any quadrangle ${\cal Q}_{ij}$ can be represented as a sum modulo 2 of
triangles and quadrangles, both of type $G(F^{n-2})$.
\end{lem}
{\bf Proof}. We span a 2-dimensional surface $S$ on the quadrangle
${\cal Q}_{ij}$ as follows. The four edges $t_i$, $t_j$, $-t_i$ and
$-t_j$ form the boundary of $S$. Recall that the vertices of the
quadrangle ${\cal Q}_{ij}$ are centers of four parallelotopes, say
$P(0)$, $P(t_i)$, $P(t_j)$ and $P(t_i+t_j)$. Hence the surface $S$
intersects a number of
parallelotopes. We choose $S$ such that it intersects boundaries of
these parallelotopes only by facets and $(n-2)$-faces, and these
intersections are transversal. Hence if the intersection $S\cap F$ is
not empty, then it is a segment or a point depending on $F$ is a facet
or an $(n-2)$-face, respectively.

These segements and points form a planar graph $\Gamma$ drawn on $S$. This
graph $\Gamma$ has 4 half-edges corresponding to the 4 facets intersected
by the 4 edges of ${\cal Q}_{ij}$. Vertices of $\Gamma$ have degrees 3 and
4 only. The dual of $\Gamma$ is a planar subgraph $G(\Gamma)$ of the graph
$G_L$. Minimal circuits of $G(\Gamma)$ are just triangles and quadrangles
of type $G(F^{n-2})$. Since $G(\Gamma)$ is planar, each edge of it
(excluding the four boundary edges) belongs to two minimal circuits. So we
obtain the wanted representation of the quadrangle ${\cal Q}_{ij}$ as a
sum of graphs $G(F^{n-2})$ modulo 2. \hfill $\Box$

\section{Pegged tilings and their duals}

McMullen in \cite{McM2} defines a pegged tiling as follows:

\vspace{2mm}
A tiling $\{Q(t):t\in{\cal T}\}$ is called {\em pegged} if with each
tile $Q(t)$, $t\in{\cal T}$, is associated a point $v^*(t)\in{\bf R}^n$,
the {\em peg} of $Q(t)$, such that if the tile $Q(t')$ is adjacent to
$Q(t)$, and so meets it on a facet $F$, then $v^*(t')-v^*(t)$ is an outer
normal vector to $Q(t)$ at the facet $F$. The equation
$x^T(v^*(t')-v^*(t))=\alpha(t,t')$ defines the hyperplane supporting $F$.

\vspace{2mm}
Note that the pegs are defined up to a shift on an arbitrary vector.
Hence we can suppose that $v^*(t^0)=0$ for some $t^0\in{\cal T}$.

Recall that parallelotopes form a tiling $\{P(t):t\in L\}$. Suppose that
this tiling is pegged. Then the peg $v^*(t)$ relates to the vertex $t$ of
the graph $G_L$. In particular, the pegs $v^*(0)=0$, $v^*(\pm t_i)$
relate to the vertices 0, $\pm t_i$, $i\in{\cal I}_P$, of the graph
$\nabla G(0)$. Since, by definition of pegs, $v^*(t_i)-v^*(0)=v^*(t_i)$
are proportional to $q_i$, and since the vectors $q_i$ are defined up to
a scalar multiple, we set $q_i=v^*(t_i)$ if the tiling $\{P(t):t\in L\}$
is pegged.

By definition of the lattice vector $t_i$, the tiles $P(t)$ and $P(t+t_i)$
are adjacent be the facet $F_i$ which is orthogonal to the facet vector
$q_i$. Hence $v^*(t+t_i)-v^*(t)=\beta_i(t)q_i$ for some scalar
$\beta_i(t)>0$, where $\beta_i(0)=1$.

In Lemma~\ref{pegq} below, we show that such defined facet vectors give
a canonical re\-pre\-sen\-ta\-ti\-on of the parallelotope $P(0)$.

\begin{lem}
\label{pegq}
The following assertions are equivalent

(i) the tiling $\{P(t);t\in L\}$ is pegged;

(ii) the parallelotope $P=P(0)$ is defined canonically.

\end{lem}
{\bf Proof}. (i)$\Rightarrow$(ii)
Let $\{t_i,t_j,t_k\}$ be a 3-belt such that $t_i=t_j+t_k$. We show that
$q_i=q_j+q_k$. Consider the hexagon of $\nabla G(0)$ corresponding to
this belt. The vertices of the hexagon are 0, $\pm t_i$, $\pm t_j$ and
$\pm t_k$. The corresponding pegs are 0, $v^*(t_r)=q_r$ and $v^*(-t_r)$,
$r\in\{i,j,k\}$. All these pegs lie in the 2-plane spanned by the facet
vectors $q_i,q_j,q_k$.

Consider the quadrangle with vertices $0=v^*(0)$, $v^*(t_j)$, $v^*(t_k)$
and $v^*(t_j+t_k)=v^*(t_i)$. The pairs of opposite edges of this quadrangle
are $(v^*(t_j)-v^*(0),v^*(t_i)-v^*(t_k))$ and
$(v^*(t_k)-v^*(0),v^*(t_i)-v^*(t_j))$. They are parallel to the vectors
$q_j$ and $q_k$, respectively. Hence this quadrangle is a parallelogram.
We have
\[q_i-q_k=v^*(t_i)-v^*(t_k)=v^*(t_j+t_k)-v^*(t_k)=
v^*(t_j)-v^*(0)=v^*(t_j)=q_j. \]
So, we obtain the wanted equality $q_i=q_j+q_k$. Since a similar reasoning
is true for every 3-belt, we see that the parallelotope $P(0)$ is defined
canonically.

(ii)$\Rightarrow$(i) Let ${\cal P}=\{t_0,t_1,...,t_s\}$ be a path
connecting the point $t=t_s$ with origin $t_0=0$. Lemma~\ref{cac} implies
that $q({\cal P})$ does not depend on the path $\cal P$, i.e.
$q({\cal P})=q(t)$.
It is easy to verify that the points $q(t)$ are pegs of the tiling
$\{P(t):t\in L\}$, i.e. $v^*(t)=q(t)$. \hfill $\Box$

\vspace{2mm}
For a pegged tiling ${\cal Q}=\{Q(t):t\in{\cal T}\}$, it is natural to
determine a tiling ${\cal Q}^*=\{Q^*(t):t\in{\cal T}^*\}$ which is
combinatorially and topologically dual to the tiling ${\cal Q}^*$ (see
\cite{McM2}, \cite{RR}).
The combinatorial duality means that, for $0\le k\le n$, there is a
one-to-one correspondence between $k$-faces of $\cal Q$ and $(n-k)$-faces
of ${\cal Q}^*$. The topological duality means that the affine spaces
spanning the corresponding $k$-face of $\cal Q$ and $(n-k)$-face of
$\cal Q^*$ are orthogonal.

So, the peg $v^*(t)$ is a vertex of the tiling $\cal Q^*$. The convex
hull of all pegs $v^*(t)$ corresponding to tiles $Q(t)\in{\cal Q}$ having
a fixed common vertex $v$ is a tile $Q^*(v)$ of $\cal Q^*$, and each tile
of $\cal Q^*$ is obtained by this way. It is proved in \cite{McM2} (see
Theorem 3.1) that $\cal Q^*$ is a tiling. Obviously, the dual tiling
$\cal Q^*$ is pegged with pegs which are vertices of the tiling $\cal Q$.
Moreover, we have $({\cal Q^*})^*=\cal Q$.

If $\cal Q$ is a pegged tiling by parallelotopes then the tiles of the
dual tiling $\cal Q^*$ are called {\em Delaunay polytopes}.
In \cite{RR}, a $k$-face of the dual tiling is called the dual convex
polytope $D^k(st)$.

Since the mutual dual tilings $\cal Q$ and $\cal Q^*$ are equivalent,
we have the following obvious assertion.
\begin{lem}
\label{pedu}
The following assertions are equivalent

(i) a tiling $\cal Q$ is pegged;

(ii) a tiling $\cal Q$ has the dual tiling $\cal Q^*$.
\end{lem}

\section{Generatrissa of a tiling}
For a tiling $\{P(t):t\in L\}$, whose tiles are parallelotopes, Voronoi
\cite{Vo} defines on the space ${\bf R}^n\otimes L$ a function $l(x;t)$.
He calls this function {\em generatrissa} and defines it as follows:

(i) $l(x;0)=0$;

(ii) if $P(t')$ is adgacent to $P(t)$ by the facet $F_i$ defined by the
equation $q_i^Tx=\alpha_i$, then
\begin{equation}
\label{rec}
l(x;t')=l(x;t)+q_i^Tx-\alpha_i.
\end{equation}

(Recall that, by (\ref{pat}), $q_i^Tx\le \alpha_i=q_i^T(t+\frac{1}{2}t_i)$
for $x\in P(t)$.)
In fact, Voronoi uses vectors $-q_i$, $i\in{\cal I}_P$, and therefore
defines $-l(x;t)$. Voronoi proves that for primitive
canonically defined parallelotopes the conditions (i) and (ii) above
determine uniquely $l(x;t)$ for each $t\in L$. The obtained generatrissa
has the following property
\begin{equation}
\label{lt0x}
l(x;t^0)\ge l(x;t) \mbox{ for all }x\in P(t^0) \mbox{ and all }t^0,t\in L,
\end{equation}
with strict inequality if $x\in{\rm int}P(t^0)$.

This property implies that the function
\[l(x)={\rm max}_{t\in L}l(x;t) \]
is a convex piecewise affine function on ${\bf R}^n$.

Consider in $(n+1)$-dimensional space ${\bf R}^n\oplus{\bf R}$ a convex
surface defined as $\{(x,z):x\in{\bf R}^n, z=l(x)\}$.
The main property of this surface is that its projection in the space
containing the tiling $\{P(t):t\in L\}$ is this tiling.

But the above definition of generatrissa works for any pegged locally
finite tiling (see, for example, \cite{McM2}).

For $t\in {\cal T}$, define $\varphi^*(t)$ as follows:

(1) $\varphi^*(t^0)=0$;

(2) $\varphi^*(t)=\sum_{i=0}^{s-1}\alpha(t^i,t^{i+1})$, where $t^s=t$,
and $P(t^0),P(t^1),...,P(t^s)=P(t)$ is a chain such that $P(t^i)$ and
$P(t^{i+1})$, $0\le i\le s-1$, are adjacent by a facet.

Here $P(t^i)$ lies but $P(t^{i+1})$ does not lie in the halfspace
$\{x\in{\bf R}^n:x^T(v^*(t^{i+1})-v^*(t^i))\le \alpha(t^i,t^{i+1})\}$.

It is proved in \cite{McM2} that the function $\varphi^*(t)$ is well defined.
In fact, it is sufficient to prove that if the chain
$P(t^0),P(t^1),...,P(t^s)=P(t^0)$ is closed, then the sum in (2) is
equal to 0. It is so if $P(t^i)$, $0\le i\le s-1$, have a common face
$F$, since, for $x_0\in F$, we have
\[\alpha(t^i,t^{i+1})=x_0^T(v^*(t^{i+1})-v^*(t^i)),\mbox{  }0\le i\le s-1.\]
Now, any closed chain can be contracted to a point avoiding $(n-3)$-faces,
and in contracting over an $(n-2)$-face we can appeal to the above reasoning.

So, the function
\begin{equation}
\label{fxt}
f(x;t)=x^Tv^*(t)-\varphi^*(t)
\end{equation}
is a generatrissa such that
\[f(x;t)>f(x;t') \mbox{ for all }x\in intP(t)
\mbox{ and all }t,t'\in{\cal T}. \]

Obviously if a tiling has a generatrissa of the form (\ref{fxt}), then
this tiling is pegged. The essential part of the papers \cite{Da} and
\cite{Au} is devoted to a proof of the following proposition.
\begin{prop}
\label{pegs}
The following assertions are equivalent for a tiling
$\{Q(t):t\in{\cal T}\}$

(i) the tiling is pegged with pegs $v^*(t)$, $t\in{\cal T}$;

(ii) the tiling has the generatrissa $f(x;t)=x^Tv^*(t)-\varphi^*(t)$.
\end{prop}

The following general result is proved in \cite{Da} and \cite{Au}.
\begin{theor}
\label{prim}
If a tiling is primitive then it is pegged and has a generatrissa.
\end{theor}
This theorem implies the main result of Voronoi \cite{Vo} asserting that
Voronoi's conjecture is true for primitive parallelotopes. The proof of
Theorem~\ref{prim} in \cite{Da} is similar to the proof of Voronoi: both
the authors construct explicitly a generatrissa. The author of \cite{Au}
constructs explicitly pegs.

\section{The case of parallelotopes}
If $P(0)$ is defined canonically, then by Lemma~\ref{pegq} the tiling
$\{P(t):t\in L\}$ is pegged with pegs $v^*(t)=q(t)$. Now, by
Proposition~\ref{pegs}, the tiling $\{P(t):t\in L\}$ by parallelotopes
has the following generatrissa
\[f(x;t)=x^Tq(t)-\varphi^*(t). \]
But for to obtain another equivalence and an explicit form of $f(x;t)$,
we use the recursion (\ref{rec}) not supposing that $P$ is defined
canonically. We want to know, when the recursion (\ref{rec}) determines
uniquely the generatrissa $l(x;t)$.

For parallelotopes, we know an explicit form of $\alpha_i$ in the
recursion (\ref{rec}). In fact, since the point $t+\frac{1}{2}t_i$
belongs to the facet $F_i$, we have $\alpha_i=q_i^T(t+\frac{1}{2}t_i)$.
Hence the recursion (\ref{rec}) takes the following form
\begin{equation}
\label{ltx}
l(x;t+t_i)=l(x;t)+q_i^T(x-(t+\frac{1}{2}t_i)).
\end{equation}
Let we have an arbitrary parallelotope $P=P(t_0)$, and let a linear
on $x$ function $l_{t^0}(x)$ be given. Then, using the recursive
expression (\ref{ltx}) and going out from $l(x;t)=l_{t^0}(x)$
by a path connecting $t^0$ and $t^1=t^0+t\in L$ we can, for all
$t^1\in L$, find a function $l(x;t^1)$ related to $L$.

Let $t=t({\cal P})=\sum_{k=1}^st_k$, where $(t_1,t_2,...,t_s)$ is an
arbitrary path $\cal P$ in $G_L$ connecting $t^0$ with $t^1$. Let
$q({\cal P})=\sum_{k=1}^sq_k$, where each $q_k$ is associated to $t_k$.
Obviously, $t({\cal P})=t^1-t^0$ does not depend on the path $\cal P$.

Using (\ref{ltx}) with $l(x;t)=l_{t^0}(x)$, and going along the path
$\cal P$ from $t^0$ we obtain
\[l(x;t^1)=l_{t^0}(x)+x^Tq({\cal P})-\phi(t^0,{\cal P}), \]
where
\begin{equation}
\label{lPx}
\phi(t^0,{\cal P})=\sum_{k=1}^sq_k^T(t^0+\sum_{r=1}^{k-1}t_r+
\frac{1}{2}t_k),
\end{equation}
and the sum $\sum_{r=1}^{k-1}t_r$ is empty for $k=1$.
(Cf., $l(x;t^1)$ with $f(x;t)$ in (\ref{fxt})).

Obviously, $q({\cal P})$ is additive on $\cal P$.
It is easy to verify that the function $\phi(t^0,{\cal P})$ is additive
on the variable $\cal P$, too. In fact, let ${\cal P}=(t_1,...,t_s)$,
${\cal P}'=(t_{s+1},...,t_w)$.
We set ${\cal P}+{\cal P}'=(t_1,...,t_s,t_{s+1},...,t_w)$. Hence the sum
of two paths is its join. Let $t=t({\cal P})$, $t'=t({\cal P}')$,
$t^1=t^0+t=t^0+\sum_{r=1}^st_r$.
Let $-{\cal P}=(-t_s,-t_{s-1},...,-t_1)$ be the path from $t^1$ to $t^0$
opposite to $\cal P$. Then
\begin{equation}
\label{add}
\phi(t^0,{\cal P}+{\cal P}')=\phi(t^0,{\cal P})+\phi(t^1,{\cal P}'),
\end{equation}
since
\[ \phi(t^0,{\cal P}+{\cal P}')=(\sum_{k=1}^s+\sum_{k=s+1}^w)q_k^T
(t^0+\sum_{r=1}^{k-1}t_r+\frac{1}{2}t_k)=\]
\[\phi(t^0,{\cal P})+\sum_{k=s+1}^wq_k^T(t^0+\sum_{r=1}^st_r+
\sum_{r=s+1}^{k-1}t_r+\frac{1}{2}t_k)=\phi(t^0,{\cal P})+
\phi(t^1,{\cal P}').\]
Now,
\[\phi(t^1,-{\cal P})=\sum_{k=1}^s(-q_{s+1-k}^T)(t^1-
\sum_{r=1}^{k-1}t_{s+1-r}-\frac{1}{2}t_{s+1-k}). \]
Since $t^1=t^0+\sum_{r=1}^st_{s+1-r}$, for $k\le s$, we have
\[t^1-\sum_{r=1}^{k-1}t_{s+1-r}=t^0+\sum_{r=1}^{s}t_{s+1-r}-
\sum_{r=1}^{k-1}t_{s+1-r}=t^0+\sum_{r=k}^st_{s+1-r}=
t^0+\sum_{r=1}^{s+1-k}t_r. \]
If we do the substitution $s+1-k\to k$ here and in $\phi(t^1,-{\cal P};x)$,
we obtain
\begin{equation}
\label{otr}
\phi(t^1,-{\cal P})=-\phi(t^0,{\cal P}).
\end{equation}

Let ${\cal C}_1={\cal P}_1+{\cal P}'$, ${\cal C}_2=-{\cal P}'+{\cal P}_2$,
where ${\cal P}_1\cap{\cal P}_2=\emptyset$. Let ${\cal P}_1$ goes from
$t^0$ to $t^1$, and the paths ${\cal P}'$, ${\cal P}_2$ go from $t^1$
to $t^0$. Then ${\cal C}_1\oplus{\cal C}_2={\cal P}_1+{\cal P}_2$ is the
sum modulo 2.
\begin{lem}
\label{C}
\[\phi(t^0,{\cal C}_1\oplus{\cal C}_2)=\phi(t^0,{\cal C}_1)+
\phi(t^0,{\cal C}_2).\]
\end{lem}
{\bf Proof}. Using (\ref{add}) and (\ref{otr}), we have
$\phi(t^0,{\cal C}_1)=\phi(t^0,{\cal P}_1+{\cal P}')=
\phi(t^0,{\cal P}_1)+\phi(t^1,{\cal P}')$ and
$\phi(t^0,{\cal C}_2)=\phi(t^0,-{\cal P}'+{\cal P}_2)=
\phi(t^0,-{\cal P}')+\phi(t^1,{\cal P}_2)=-\phi(t^1,{\cal P}')+
\phi(t^1,{\cal P}_2)$. Hence
\[\phi(t^0,{\cal C}_1)+\phi(t^0,{\cal C}_2)=
\phi(t^0,{\cal P}_1)+\phi(t^1,{\cal P}_2)=
\phi(t^0,{\cal P}_1+{\cal P}_2)=\phi(t^0,{\cal C}_1\oplus{\cal C}_2).
\hfill \Box \]

Using Lemmas~\ref{C} and \ref{quad}, for any quadrangle ${\cal Q}_{ij}$,
we can represent $\phi(t^0,{\cal Q}_{ij})$ as a sum of functions
$\phi(t,G(F^{n-2}))$.

\begin{lem}
\label{val}
Let $\phi(t^0,{\cal P})$ be given by (\ref{lPx}). Then, for
$i,j\in{\cal I}_P$ and ${\cal P}={\cal Q}_{ij}$,
\[\phi(t^0,{\cal Q}_{ij})=q_j^Tt_i-q_i^Tt_j. \]
\end{lem}
{\bf Proof}. Using the equality (\ref{lPx}) for
${\cal P}={\cal Q}_{ij}=(t_i,t_j,-t_i,-t_j)$, we obtain
\[\phi(t^0,{\cal Q}_{ij})=q_i^T(t^0+\frac{1}{2}t_i)+
q_j^T(t^0+t_i+\frac{1}{2}t_j)-\]
\[-q_i^T(t^0+t_i+t_j-\frac{1}{2}t_i)-q_j^T(t^0+t_i+t_j-t_i-
\frac{1}{2}t_j)=q_j^Tt_i-q_i^Tt_j. \hfill \Box \]

Lemma~\ref{val} implies the following result.

\begin{lem}
\label{Qij}
Let $i,j\in{\cal I}_P$, and let $\phi(t^0,{\cal P})$ be given by (\ref{lPx}).
The following assertions are equivalent:

(i) $\phi(t^0,{\cal Q}_{ij})=0$;

(ii) $q_i^Tt_j=q_j^Tt_i$. \hfill $\Box$
\end{lem}

Let ${\cal P}=(t_1,t_2,...,t_s)$ be a path. Recall that
$t({\cal P})=\sum_{k=1}^st_k$, $q({\cal P})=\sum_{k=1}^sq_k$, and
${\cal I}({\cal P})$ is the set of all indices $i\in{\cal I}_P$ of $t_i$
in the path $\cal P$.
\begin{lem}
\label{tPx}
Let the equalities $q_i^Tt_j=q_j^Tt_i$ hold for all pairs
$i,j\in{\cal I}({\cal P})$. Then the function $\phi(t^0,{\cal P})$ is given
by the expression
\begin{equation}
\label{lP}
\phi(t^0,{\cal P})=q^T({\cal P})(t^0+\frac{1}{2}t({\cal P})).
\end{equation}
\end{lem}
{\bf Proof}.
Using the equalities $q_i^Tt_j=t_i^Tq_j$ for $i,j\in{\cal I}({\cal P})$,
we obtain the equality
\[q_k^T\sum_{r=1}^{k-1}t_r=t_k^T\sum_{r=1}^{k-1}q_r. \]
Using this equality and setting $\sum_{k=1}^sq_k=q({\cal P})$, we rewrite
$\phi(t^0,{\cal P})$ from (\ref{lPx}) as follows
\[\phi(t^0,{\cal P})=\sum_{k=1}^sq_k^Tt^0+
\frac{1}{2}\sum_{k=1}^sq_k^T\sum_{r=1}^{k-1}t_r+
\frac{1}{2}\sum_{k=1}^sq_k^T(\sum_{r=1}^{k-1}t_r+t_k)= \]
\[=q^T({\cal P})t^0+\frac{1}{2}\sum_{k=1}^st_k^T\sum_{r=1}^{k-1}q_r+
\frac{1}{2}\sum_{k=1}^sq_k^T\sum_{r=1}^{k}t_r. \]
A permutation of the summation in the first double sum, gives
\[\frac{1}{2}\sum_{k=1}^st_k^T\sum_{r=1}^{k-1}q_r=
\frac{1}{2}\sum_{k=1}^{s-1}q_k^T\sum_{r=k+1}^{s}t_r. \]
Recall that $\sum_{k=1}^st_k=t({\cal P})$. Hence $\phi(t^0,{\cal P})$
takes the form
\[\phi(t^0,{\cal P})=q^T({\cal P})t^0+\frac{1}{2}\sum_{k=1}^sq_k^T
\sum_{r=1}^{s}t_r=q^T({\cal P})(t^0+\frac{1}{2}t({\cal P})). \]
So, we obtain the wanted expression. \hfill $\Box$

\vspace{3mm}
Using Lemma~\ref{tPx}, we can prove the following important result.
\begin{lem}
\label{can}
The following assertions are equivalent:

(i) a parallelotope $P$ is defined canonically;

(ii) $q_i^Tt_j=q_j^Tt_i$ for all $i,j\in{\cal I}_P$.
\end{lem}
{\bf Proof}. (i)$\Rightarrow$(ii). If $P$ is defined canonically, then
Proposition~\ref{belt} assert that $q_i^Tt_j=q_j^Tt_i$
for $i,j\in{\cal I}(G(F^{n-2}))$. Obviously, $q({\cal Q}_{ij})=0$, and
by Lemma~\ref{tri}, $q(\triangle)=0$ if $\triangle=G(F^{n-2})$. Hence
$q(G(F^{n-2}))=0$ and Lemma~\ref{tPx} implies that $\phi(t^0,G(F^{n-2}))=0$.
Lemma~\ref{quad} asserts that each quadrangle ${\cal Q}_{ij}$ is a sum
modulo 2 of circuits $G(F^{n-2})$. By Lemma~\ref{C} we obtain that
$\phi(t^0,{\cal Q}_{ij})=0$. Now, Lemma~\ref{Qij} gives the
wanted equality for all pairs $i,j$.

(ii)$\Rightarrow$(i) By Lemma~\ref{qt}, $q_i=Dt_i$ for all $i\in{\cal I}_P$.
Hence any linear equality between lattice vectors implies the corresponding
equality between the associated facet vectors. This implies that the facet
vectors are defined canonically with respect to each belt of $P$, i.e.,
$P$ is defined canonically. \hfill $\Box$

The obtained function $\phi(t^0,{\cal P})$ depends on the chosen path
$\cal P$. Since $\phi(t^0,{\cal P})$ satisfies the condition (\ref{otr}),
it does not depend on $\cal P$ if and only if $\phi(t^0,{\cal C})=0$ for
every circuit $\cal C$. In Lemma~\ref{lin} below we give conditions when
this property is true.
\begin{lem}
\label{lin}
Let $\cal C$ be a circuit and $t^0\in\cal C$. The following
assertions are equivalent:

(i) for all $t^0$ and all circuits ${\cal C}\ni t^0$, the function
$\phi(t^0,{\cal C})=0$;

(ii) the equalities $q_i^Tt_j=q_j^Tt_i$ hold for all pairs
$i,j\in{\cal I}_P$.
\end{lem}
{\bf Proof}. (i)$\Rightarrow$(ii) The item (i) implies that
$\phi(t^0,{\cal Q}_{ij})=0$ for any quadrangle ${\cal Q}_{ij}$. Now
Lemmas~\ref{Qij} gives wanted implication.

(ii)$\Rightarrow$(i) Recall that if $q^T_it_j=q^T_jt_i$
for all $i,j\in {\cal I}_P$, then $q_i=Dt_i$. Since the sum
$\sum_{k=1}^st_k=t({\cal C})=0$ for any circuit $\cal C$, the associated
vector $q({\cal C})=\sum_{k=1}^sq_k=\sum_{k=1}^sDt_k=Dt({\cal C})$
is also equal to zero. Hence
$\phi(t^0,{\cal C})=q^T({\cal C})(t^0+\frac{1}{2}t({\cal C}))=0$.
\hfill $\Box$

\section{Generatrissa for parallelotopes}

By Lemma~\ref{lin} the generatrissa $l(x;t)=q^T(t)(x-\frac{1}{2}t)$
is uniquely determined for all $x\in{\bf R}^n$ and $t\in L$.
This is the function $q^T({\cal P})x-\phi(0,{\cal P})$, where
$q({\cal P})=q(t)$ and $\phi(0,{\cal P})=\frac{1}{2}q^T(t)t$, both do
not depend on $\cal P$. In other words, we have the following assertion
\begin{lem}
\label{gen}
The following assertions are equivalent:

(i) the function $q^T(t)(x-\frac{1}{2}t)$, $t\in L$, is the generatrissa
of the tiling obtained by translations of a parallelotope $P$;

(ii) the equalities $q_i^Tt_j=q_j^Tt_i$ hold for all $i,j\in{\cal I}_P$.
\hfill $\Box$

\end{lem}

Recall that the facet vector $q(t)$ is associated to the lattice vector
$t$, when (\ref{QT}) is true. Hence, by Lemma~\ref{qt}, $q(t)=Dt$ and we
see that $l(x;t)$ is given by
\[l(x;t)=t^TD(x-\frac{1}{2}t). \]

We see that $l(x;0)=0$. Using standard arguments of \cite{Vo} and
\cite{Zh}, one can prove that $l(x;t)\le 0$ for all $x\in P(0)$ and all
$t\in L$, with strict inequality for $x\in{\rm int}P(0)$. Hence
(\ref{lt0x}) (with $t^0=0$) is true. But Lemma~\ref{qt} does not asserts
that the matrix $D$ is positive definite. The next lemma proves that $D$
is positive definite if the vectors $q_i$ and $t_i$ are related to a
parallelotope.
\begin{lem}
\label{pos}
Let $P(0)$ be a parallelotope with center in origin, and let (\ref{QT})
is true. Let $D$ be a non-singular matrix. The following assertions are
equivalent:

(i) the inequality $l(x;t)=t^TD(x-\frac{1}{2}t)< 0$ holds for all
$x\in {\rm int}P(0)$ and all $t\in L$, $t\not=0$;

(ii) the matrix $D$ is positive definite.
\end{lem}
{\bf Proof}. (i)$\Rightarrow$(ii) For $x=0$ and $t\not=0$, we have
$-l(0;t)=\frac{1}{2}t^TDt> 0$. This inequality holds for all vectors
$t\in L-\{0\}$. Any rational combination of basic vectors of $L$ is equal to
$\frac{1}{p}t$ for some $t\in L$ and an integer $p$. We obtain that the
above inequality holds also for rational vectors. By continuity, this
inequality holds for all $t\in {\bf R}^n$. This implies, that the matrix
$D$ is positive definite.

(ii)$\Rightarrow$(i) Let $D$ be positive definite. Consider the quadratic
function $x^TDx-2l(x;t)=(x-t)^TD(x-t)=f(x-t)$. Using the quadratic
form $f(x)$ we can define the parallelotope
\[P_f(t^0)=\{x\in{\bf R}^n:l(x;t^0)-l(x;t)\ge 0, \mbox{  }t\in L\}. \]
We show that $P_f(t^0)=P(t^0)$. Then for $t^0=0$ we will have
$l(x;t)\le 0$, $t\in L-\{0\}$, and $l(x;t)<0$ for interior points of
$P(t^0)$. The infinite system inequalities
\[l(x;t^0)-l(x;t)\ge 0, \mbox{  }t\in L, \]
contains the subsystem for $t=t^0\pm t_i$, $i\in {\cal I}_P$, describing
the parallelotope $P(t^0)$ of type (\ref{pat}). Hence 
$P_f(t^0) \subseteq P(t^0)$ for every $t^0\in L$. But we have here an
equality. In fact, if there is $t\in L$ such that $P_f(t)\subset P(t)$
strictly, then there is an adjacent to $P_f(t)$ parallelotope $P_f(t')$
such that $P_f(t')$ and $P(t)$ have a common interior point. Since $P(t')$
contains $P_f(t')$, the parallelotopes $P(t)$ and $P(t')$ have a common
interior point. This is a contradiction. \hfill $\Box$

The above proof of Lemma~\ref{pos} is also a proof of the following
\begin{lem}
\label{fD}
The following assertions are equivalent

(i) the parallelotope $P$ is a Voronoi polytope with respect to the
positive quadratic form $f(x)=x^TDx$;

(ii) the function $l(x;t)=t^TD(x-\frac{1}{2}t)$ is a generatrissa of
the tiling $\{P(t):t\in L\}$.
\end{lem}

\section{Main Theorem}
If we collect Lemmas~\ref{ADq}, \ref{qt}, \ref{pegq}, \ref{pedu},
Proposition~\ref{pegs} and Lemmas~\ref{can}, \ref{gen},
\ref{fD} together, we obtain
\begin{theor}
\label{main}
Let $P$ be a parallelotope defined by facet vectors $q_i$, and defining
lattice vectors $t_i$, $i\in{\cal I}_P$. Let $A$ be a non-degenerate
$n\times n$ matrix, and $D=A^TA$.
Then the following assertions are equivalent:

(i) Voronoi's conjecture holds for $P$, i.e. the affine transformation
$x\to Ax$ transforms the parallelotope $P$ into a Voronoi polytope;

(ii) the parallelotope $P$ is a Voronoi polytope with respect to the
positive quadratic form $f(x)=x^TDx$;

(iii) the parallelotope $P$ is defined canonically;

(iv) the equality $q_i=Dt_i$ holds for all $i\in {\cal I}_P$;

(v) for all pairs $i,j\in {\cal I}_P$, the equalities $t_i^Tq_j=q_i^Tt_j$
hold;

(vi) the tiling $\{P(t):t\in L\}$ is pegged with pegs $v^*(t)=Dt$;

(vii) the function $l(x;t)=t^TD(x-\frac{1}{2}t)$ is a generatrissa of
the tiling $\{P(t):t\in L\}$;

(viii) the tiling $\{P(t):t\in L\}$ has a dual tiling.
\end{theor}

Theorem~\ref{main} implies a result of \cite{MRS} that there is a unique
(up to isomorphism of $P$) map which transforms a primitive parallelorope
$P$ into a Voronoi polytope. Our Theorem~\ref{main} gives an explicit
matrix $A$ of the corresponding affine map. Therefore we have the following
\begin{prop}
If a parallelotope $P$ is affinely equivalent to a Voronoi polytope, then
this affinity is uniquely (up to the aftomorphism of $P$) determined by
the parallelotope $P$.
\end{prop}

M.Deza \hfill V.Grishukhin, \\
Ecole Normale Sup\'erieure \hfill CEMI Russian Academy of Sciences\\
r.Ulm 45 \hfill Nakhimovskii prosp.47\\
Paris, France \hfill Moscow, Russia \\
e-mail: deza@ens.fr \hfill e-mail: grishuhn@cemi.rssi.ru\\

\end{document}